\documentclass{amsart}
\usepackage{amssymb}
\usepackage{amsfonts}

\newtheorem{theorem}{Theorem}
\theoremstyle{plain}

\newtheorem{corollary}{Corollary}

\newtheorem{lemma}{Lemma}

\newtheorem{remark}{Remark}

\numberwithin{equation}{section}
\input{tcilatex}

\begin{document}
\title[\v{C}eby\v{s}ev Functional]{Bounding the \v{C}eby\v{s}ev Functional
for a Pair of Sequences in Inner Product Spaces}
\author{S.S. Dragomir}
\address{School of Computer Science and Mathematics\\
Victoria University of Technology\\
PO Box 14428, MCMC 8001\\
Victoria, Australia.}
\email{sever.dragomir@vu.edu.au}
\urladdr{http://rgmia.vu.edu.au/SSDragomirWeb.html}
\date{September 16, 2003.}

\begin{abstract}
Some new bounds for the \v{C}eby\v{s}ev functional of a pair of vectors in
inner product spaces are pointed out. Reverses for the celebrated Jensen's
inequality for convex functions defined on inner product spaces are given as
well.
\end{abstract}

\maketitle

\section{Introduction}

Let $\left( H;\left\langle \cdot ,\cdot \right\rangle \right) $ be an inner
product over the real or complex number field $\mathbb{K}$. For $\mathbf{p}%
=\left( p_{1},\dots ,p_{n}\right) \in \mathbb{R}^{n}$ and $\mathbf{x}=\left(
x_{1},\dots ,x_{n}\right) ,$ $\mathbf{y}=\left( y_{1},\dots ,y_{n}\right)
\in H^{n},$ define the \textit{\v{C}eby\v{s}ev functional}%
\begin{equation}
T_{n}\left( \mathbf{p};\mathbf{x},\mathbf{y}\right)
:=P_{n}\sum_{i=1}^{n}p_{i}\left\langle x_{i},y_{i}\right\rangle
-\left\langle \sum_{i=1}^{n}p_{i}x_{i},\sum_{i=1}^{n}p_{i}y_{i}\right\rangle
,  \label{1.1}
\end{equation}%
where $P_{n}:=\sum_{i=1}^{n}p_{i}.$

The following Gr\"{u}ss type inequality has been obtained in \cite{SSD1}.

\begin{theorem}
\label{t1.1}Let $H,$ $\mathbf{x},\mathbf{y}$ be as above and $p_{i}\geq 0$ \ 
$\left( i\in \left\{ 1,\dots ,n\right\} \right) $ with $%
\sum_{i=1}^{n}p_{i}=1,$ i.e., $\mathbf{p}$ is a probability sequence. If $%
x,X,y,Y\in H$ are such that%
\begin{equation}
\func{Re}\left\langle X-x_{i},x_{i}-x\right\rangle \geq 0,\ \ \ \func{Re}%
\left\langle Y-y_{i},y_{i}-y\right\rangle \geq 0  \label{1.2}
\end{equation}%
for each $i\in \left\{ 1,\dots ,n\right\} ,$ or equivalently, (see \cite%
{SSD2})%
\begin{equation}
\left\Vert x_{i}-\frac{x+X}{2}\right\Vert \leq \frac{1}{2}\left\Vert
X-x\right\Vert ,\ \ \ \ \ \ \left\Vert y_{i}-\frac{y+Y}{2}\right\Vert \leq 
\frac{1}{2}\left\Vert Y-y\right\Vert  \label{1.2a}
\end{equation}%
for each $i\in \left\{ 1,\dots ,n\right\} ,$ then we have the inequality%
\begin{equation}
\left\vert T_{n}\left( \mathbf{p};\mathbf{x},\mathbf{y}\right) \right\vert
\leq \frac{1}{4}\left\Vert X-x\right\Vert \left\Vert Y-y\right\Vert .
\label{1.3}
\end{equation}%
The constant $\frac{1}{4}$ is best possible in the sense that it cannot be
replaced by a smaller constant.
\end{theorem}

In \cite{SSD3}, the following Gr\"{u}ss type inequality for the forward
difference of vectors was established.

\begin{theorem}
\label{t1.2}Let $\mathbf{x}=\left( x_{1},\dots ,x_{n}\right) ,$ $\mathbf{y}%
=\left( y_{1},\dots ,y_{n}\right) \in H^{n}$ and $\mathbf{p}\in \mathbb{R}%
_{+}^{n}$ be a probability sequence. Then one has the inequality:%
\begin{align}
& \left\vert T_{n}\left( \mathbf{p};\mathbf{x},\mathbf{y}\right) \right\vert
\label{1.4} \\
& \leq \left\{ 
\begin{array}{l}
\left[ \sum\limits_{i=1}^{n}i^{2}p_{i}-\left(
\sum\limits_{i=1}^{n}ip_{i}\right) ^{2}\right] \max\limits_{1\leq k\leq
n-1}\left\Vert \Delta x_{k}\right\Vert \max\limits_{1\leq k\leq
n-1}\left\Vert \Delta y_{k}\right\Vert ; \\ 
\\ 
\sum\limits_{1\leq j<i\leq n}p_{i}p_{j}\left( i-j\right) \left(
\sum\limits_{k=1}^{n-1}\left\Vert \Delta x_{k}\right\Vert ^{p}\right) ^{%
\frac{1}{p}}\left( \sum\limits_{k=1}^{n-1}\left\Vert \Delta y_{k}\right\Vert
^{q}\right) ^{\frac{1}{q}} \\ 
\hfill \text{ \ if \ }p>1,\ \frac{1}{p}+\frac{1}{q}+1 \\ 
\\ 
\dfrac{1}{2}\left[ \sum\limits_{i=1}^{n}p_{i}\left( 1-p_{i}\right) \right]
\sum\limits_{k=1}^{n-1}\left\Vert \Delta x_{k}\right\Vert
\sum\limits_{k=1}^{n-1}\left\Vert \Delta y_{k}\right\Vert .%
\end{array}%
\right.  \notag
\end{align}%
The constants $1,1$ and $\frac{1}{2}$ in the right hand side of inequality (%
\ref{1.4}) are best in the sense that they cannot be replaced by smaller
constants.
\end{theorem}

Another result is incorporated in the following theorem (see \cite{SSD2}).

\begin{theorem}
\label{t1.3}Let $\mathbf{x},\mathbf{y}$ and $\mathbf{p}$ be as in Theorem %
\ref{t1.2}. If there exist $x,X\in H$ such that%
\begin{equation}
\func{Re}\left\langle X-x_{i},x_{i}-x\right\rangle \geq 0\ \ \text{ for each
\ }i\in \left\{ 1,\dots ,n\right\} ,  \label{1.5}
\end{equation}%
or, equivalently,%
\begin{equation}
\left\Vert x_{i}-\frac{x+X}{2}\right\Vert \leq \frac{1}{2}\left\Vert
X-x\right\Vert \ \ \text{ for each \ }i\in \left\{ 1,\dots ,n\right\} ,
\label{1.6}
\end{equation}%
then one has the inequality%
\begin{align}
\left\vert T_{n}\left( \mathbf{p};\mathbf{x},\mathbf{y}\right) \right\vert &
\leq \frac{1}{2}\left\Vert X-x\right\Vert
\sum\limits_{i=1}^{n}p_{i}\left\Vert
y_{i}-\sum_{j=1}^{n}p_{j}y_{j}\right\Vert  \label{1.7} \\
& \leq \frac{1}{2}\left\Vert X-x\right\Vert \left[ \sum_{i=1}^{n}p_{i}\left%
\Vert y_{i}\right\Vert ^{2}-\left\Vert
\sum\limits_{i=1}^{n}p_{i}y_{i}\right\Vert ^{2}\right] ^{\frac{1}{2}}. 
\notag
\end{align}%
The constant $\frac{1}{2}$ is best possible in the first and second
inequalities in the sense that it cannot be replaced by a smaller constant.
\end{theorem}

\begin{remark}
If $\mathbf{x}$ and $\mathbf{y}$ satisfy the assumptions of Theorem \ref%
{t1.1}, then we have the following sequence of inequalities improving the Gr%
\"{u}ss inequality (\ref{1.3}):%
\begin{align}
\left\vert T_{n}\left( \mathbf{p};\mathbf{x},\mathbf{y}\right) \right\vert &
\leq \frac{1}{2}\left\Vert X-x\right\Vert \sum_{i=1}^{n}p_{i}\left\Vert
y_{i}-\sum_{j=1}^{n}p_{j}y_{j}\right\Vert  \label{1.8} \\
& \leq \frac{1}{2}\left\Vert X-x\right\Vert \left(
\sum_{i=1}^{n}p_{i}\left\Vert y_{i}\right\Vert ^{2}-\left\Vert
\sum\limits_{i=1}^{n}p_{i}y_{i}\right\Vert ^{2}\right) ^{\frac{1}{2}}  \notag
\\
& \leq \frac{1}{4}\left\Vert X-x\right\Vert \left\Vert Y-y\right\Vert . 
\notag
\end{align}
\end{remark}

Now, if we consider the \v{C}eby\v{s}ev functional for the uniform
probability distribution $u=\left( \frac{1}{n},\dots ,\frac{1}{n}\right) ,$%
\begin{equation*}
T_{n}\left( \mathbf{x},\mathbf{y}\right) :=\frac{1}{n}\sum_{i=1}^{n}\left%
\langle x_{i},y_{i}\right\rangle -\left\langle \frac{1}{n}%
\sum_{i=1}^{n}x_{i},\frac{1}{n}\sum_{i=1}^{n}y_{i}\right\rangle ,
\end{equation*}%
then, with the assumptions of Theorem \ref{t1.1}, we have%
\begin{equation}
\left\vert T_{n}\left( \mathbf{x},\mathbf{y}\right) \right\vert \leq \frac{1%
}{4}\left\Vert X-x\right\Vert \left\Vert Y-y\right\Vert .  \label{1.9}
\end{equation}

Theorem \ref{t1.2} will provide the following inequalities%
\begin{multline}
\left\vert T_{n}\left( \mathbf{x},\mathbf{y}\right) \right\vert  \label{1.10}
\\
\leq \left\{ 
\begin{array}{l}
\dfrac{1}{12}\left( n^{2}-1\right) \max\limits_{1\leq k\leq n-1}\left\Vert
\Delta x_{k}\right\Vert \max\limits_{1\leq k\leq n-1}\left\Vert \Delta
y_{k}\right\Vert ; \\ 
\\ 
\dfrac{1}{6}\left( n-\dfrac{1}{n}\right) \left(
\sum\limits_{k=1}^{n-1}\left\Vert \Delta x_{k}\right\Vert ^{p}\right) ^{%
\frac{1}{p}}\left( \sum\limits_{k=1}^{n-1}\left\Vert \Delta y_{k}\right\Vert
^{q}\right) ^{\frac{1}{q}}\hfill \text{ \ if \ }p>1,\ \frac{1}{p}+\frac{1}{q}%
+1; \\ 
\\ 
\dfrac{1}{2}\left( 1-\dfrac{1}{n}\right) \sum\limits_{k=1}^{n-1}\left\Vert
\Delta x_{k}\right\Vert \sum\limits_{k=1}^{n-1}\left\Vert \Delta
y_{k}\right\Vert .%
\end{array}%
\right.
\end{multline}%
Here the constants $\frac{1}{12},$ $\frac{1}{6}$ and $\frac{1}{2}$ are best
possible in the above sense.

Finally, from (\ref{1.8}), we have%
\begin{align}
\left\vert T_{n}\left( \mathbf{x},\mathbf{y}\right) \right\vert & \leq \frac{%
1}{2n}\left\Vert X-x\right\Vert \sum_{i=1}^{n}\left\Vert y_{i}-\dfrac{1}{n}%
\sum_{j=1}^{n}y_{j}\right\Vert  \label{1.11} \\
& \leq \frac{1}{2}\left\Vert X-x\right\Vert \left( \frac{1}{n}%
\sum_{i=1}^{n}\left\Vert y_{i}\right\Vert ^{2}-\left\Vert \frac{1}{n}%
\sum_{i=1}^{n}y_{i}\right\Vert ^{2}\right) ^{\frac{1}{2}}  \notag \\
& \leq \frac{1}{4}\left\Vert X-x\right\Vert \left\Vert Y-y\right\Vert . 
\notag
\end{align}

It is the main aim of this paper to point out other bounds for the \v{C}eby%
\v{s}ev functionals $T_{n}\left( \mathbf{p};\mathbf{x},\mathbf{y}\right) $
and $T_{n}\left( \mathbf{x},\mathbf{y}\right) .$ Applications for Jensen's
inequality for convex functions defined on inner product spaces are given as
well.

\section{Identities for Inner Products}

For $\mathbf{p}=\left( p_{1},\dots ,p_{n}\right) \in \mathbb{R}^{n}$ and $%
\mathbf{a}=\left( a_{1},\dots ,a_{n}\right) \in H^{n}$ we define%
\begin{equation*}
P_{i}:=\sum_{k=1}^{i}p_{k},\ \ \ \ \ \ \bar{P}_{i}=P_{n}-P_{i},\ \ \ \ \
i\in \left\{ 1,\dots ,n-1\right\}
\end{equation*}%
and the vectors%
\begin{equation*}
A_{i}\left( \mathbf{p}\right) =\sum_{k=1}^{i}p_{k}a_{k},\ \ \ \ \ \ \bar{A}%
_{i}\left( \mathbf{p}\right) =A_{n}\left( \mathbf{p}\right) -A_{i}\left( 
\mathbf{p}\right)
\end{equation*}%
for $i\in \left\{ 1,\dots ,n-1\right\} .$

The following result holds.

\begin{theorem}
\label{t2.1}Let $\left( H;\left\langle \cdot ,\cdot \right\rangle \right) $
be an inner product space over $\mathbb{K}$, $\mathbf{p}=\left( p_{1},\dots
,p_{n}\right) \in \mathbb{R}^{n}$ and $\mathbf{a}=\left( a_{1},\dots
,a_{n}\right) ,\mathbf{b}=\left( b_{1},\dots ,b_{n}\right) \in H^{n}.$ Then
we have the identities%
\begin{align}
T_{n}\left( \mathbf{p};\mathbf{a},\mathbf{b}\right) &
=\sum_{i=1}^{n-1}\left\langle P_{i}A_{n}\left( \mathbf{p}\right)
-P_{n}A_{i}\left( \mathbf{p}\right) ,\Delta b_{i}\right\rangle  \label{2.1}
\\
& =P_{n}\sum_{i=1}^{n-1}P_{i}\left\langle \frac{1}{P_{n}}A_{n}\left( \mathbf{%
p}\right) -\frac{1}{P_{i}}A_{i}\left( \mathbf{p}\right) ,\Delta
b_{i}\right\rangle  \notag \\
\text{(if }P_{i}& \neq 0,\ i\in \left\{ 1,\dots ,n\right\} \text{)}  \notag
\\
& =\sum_{i=1}^{n-1}P_{i}\bar{P}_{i}\left\langle \frac{1}{\bar{P}_{i}}\bar{A}%
_{i}\left( \mathbf{p}\right) -\frac{1}{P_{i}}A_{i}\left( \mathbf{p}\right)
,\Delta b_{i}\right\rangle  \notag \\
\text{(if }P_{i},\bar{P}_{i}& \neq 0,\ i\in \left\{ 1,\dots ,n-1\right\} 
\text{),}  \notag
\end{align}%
where $\Delta x_{i}=x_{i+1}-x_{i}$ $\left( i\in \left\{ 1,\dots ,n-1\right\}
\right) $ is the forward difference.
\end{theorem}

\begin{proof}
We use the following summation by parts formula for vectors in inner product
spaces%
\begin{equation}
\sum_{l=p}^{q-1}\left\langle d_{l},\Delta v_{l}\right\rangle =\left\langle
d_{l},v_{l}\right\rangle \big|_{p}^{q}-\sum_{l=p}^{q-1}\left\langle
v_{l+1},\Delta d_{l}\right\rangle ,  \label{2.2}
\end{equation}%
where $d_{l},$ $v_{l}$ are vectors in $H,$ $l=p,\dots ,q$ ($q>p;$ $p,q$ are
natural numbers).

If we choose in (\ref{2.2}), $p=1,$ $q=n,$ $d_{i}=P_{i}A_{n}\left( \mathbf{p}%
\right) -P_{n}A_{i}\left( \mathbf{p}\right) $ and $v_{i}=b_{i}$ $\left( i\in
\left\{ 1,\dots ,n-1\right\} \right) ,$ then we get%
\begin{eqnarray*}
&&\sum_{i=1}^{n-1}\left\langle P_{i}A_{n}\left( \mathbf{p}\right)
-P_{n}A_{i}\left( \mathbf{p}\right) ,\Delta b_{i}\right\rangle \\
&=&\left\langle P_{i}A_{n}\left( \mathbf{p}\right) -P_{n}A_{i}\left( \mathbf{%
p}\right) ,b_{i}\right\rangle \big|_{1}^{n}-\sum_{i=1}^{n-1}\left\langle
\Delta \left( P_{i}A_{n}\left( \mathbf{p}\right) -P_{n}A_{i}\left( \mathbf{p}%
\right) \right) ,b_{i+1}\right\rangle \\
&=&\left\langle P_{n}A_{n}\left( \mathbf{p}\right) -P_{n}A_{n}\left( \mathbf{%
p}\right) ,b_{n}\right\rangle -\left\langle P_{1}A_{n}\left( \mathbf{p}%
\right) -P_{n}A_{1}\left( \mathbf{p}\right) ,b_{1}\right\rangle \\
&&-\sum_{i=1}^{n-1}\left\langle P_{i+1}A_{n}\left( \mathbf{p}\right)
-P_{n}A_{i+1}\left( \mathbf{p}\right) -P_{i}A_{n}\left( \mathbf{p}\right)
+P_{n}A_{i}\left( \mathbf{p}\right) ,b_{i+1}\right\rangle \\
&=&P_{n}p_{1}\left\langle a_{1},x_{1}\right\rangle -p_{1}\left\langle
A_{n}\left( \mathbf{p}\right) ,b_{1}\right\rangle -\left\langle A_{n}\left( 
\mathbf{p}\right) ,\sum_{i=1}^{n-1}p_{i+1}b_{i+1}\right\rangle \\
&&+P_{n}\sum_{i=1}^{n-1}p_{i+1}\left\langle a_{i+1},b_{i+1}\right\rangle \\
&=&P_{n}\sum_{i=1}^{n}p_{i}\left\langle a_{i},b_{i}\right\rangle
-\left\langle \sum_{i=1}^{n}p_{i}a_{i},\sum_{i=1}^{n}p_{i}b_{i}\right\rangle
\\
&=&T_{n}\left( \mathbf{p};\mathbf{a},\mathbf{b}\right) ,
\end{eqnarray*}%
proving the first identity in (\ref{2.1}).

The second and third identities are obvious and we omit the details.
\end{proof}

The following lemma is of interest in itself.

\begin{lemma}
\label{l2.2}Let $\mathbf{p}=\left( p_{1},\dots ,p_{n}\right) \in \mathbb{R}%
^{n}$ and $\mathbf{a}=\left( a_{1},\dots ,a_{n}\right) \in H.$ Then we have
the equality%
\begin{equation}
P_{i}A_{n}\left( \mathbf{p}\right) -P_{n}A_{i}\left( \mathbf{p}\right)
=\sum_{j=1}^{n-1}P_{\min \left\{ i,j\right\} }\bar{P}_{\max \left\{
i,j\right\} }\Delta a_{j}  \label{2.3}
\end{equation}%
for each $i\in \left\{ 1,\dots ,n-1\right\} .$
\end{lemma}

\begin{proof}
Define, for $i\in \left\{ 1,\dots ,n-1\right\} ,$ the vector%
\begin{equation*}
K\left( i\right) :=\sum_{j=1}^{n-1}P_{\min \left\{ i,j\right\} }\bar{P}%
_{\max \left\{ i,j\right\} }\cdot \Delta a_{j}.
\end{equation*}%
We have%
\begin{align}
K\left( i\right) & =\sum_{j=1}^{i}P_{\min \left\{ i,j\right\} }\bar{P}_{\max
\left\{ i,j\right\} }\cdot \Delta a_{j}+\sum_{j=i+1}^{n-1}P_{\min \left\{
i,j\right\} }\bar{P}_{\max \left\{ i,j\right\} }\cdot \Delta a_{j}
\label{2.4} \\
& =\sum_{j=1}^{i}P_{j}\bar{P}_{i}\cdot \Delta a_{j}+\sum_{j=i+1}^{n-1}P_{i}%
\bar{P}_{j}\cdot \Delta a_{j}  \notag \\
& =\bar{P}_{i}\sum_{j=1}^{i}P_{j}\cdot \Delta a_{j}+P_{i}\sum_{j=i+1}^{n-1}%
\bar{P}_{j}\cdot \Delta a_{j}.  \notag
\end{align}%
Using the summation by parts formula, we have%
\begin{align}
\sum_{j=1}^{i}P_{j}\cdot \Delta a_{j}& =P_{j}a_{j}\big|_{1}^{i+1}-%
\sum_{j=1}^{i}\left( P_{j+1}-P_{j}\right) a_{j+1}  \label{2.5} \\
& =P_{i+1}a_{i+1}-p_{1}a_{1}-\sum_{j=1}^{i}p_{j+1}a_{j+1}  \notag \\
& =P_{i+1}a_{i+1}-\sum_{j=1}^{i+1}p_{j}a_{j}  \notag
\end{align}%
and%
\begin{align}
\sum_{j=i+1}^{n-1}\bar{P}_{j}\cdot \Delta a_{j}& =\bar{P}_{j}a_{j}\big|%
_{i+1}^{n}-\sum_{j=i+1}^{n-1}\left( \bar{P}_{j+1}-\bar{P}_{j}\right) a_{j+1}
\label{2.6} \\
& =\bar{P}_{n}a_{n}-\bar{P}_{i+1}a_{i+1}-\sum_{j=i+1}^{n-1}\left(
P_{n}-P_{j+1}-P_{n}+P_{j}\right) a_{j+1}  \notag \\
& =-\bar{P}_{i+1}a_{i+1}+\sum_{j=i+1}^{n-1}p_{j+1}a_{j+1}.  \notag
\end{align}%
Using (\ref{2.5}) and (\ref{2.6}), we have%
\begin{align*}
K\left( i\right) & =\bar{P}_{i}\left(
P_{i+1}a_{i+1}-\sum_{j=1}^{i+1}p_{j}a_{j}\right) +P_{i}\left(
\sum_{j=i+1}^{n-1}p_{j+1}a_{j+1}-\bar{P}_{i+1}a_{i+1}\right) \\
& =\bar{P}_{i}P_{i+1}a_{i+1}-\bar{P}_{i}\bar{P}_{i+1}a_{i+1}-\bar{P}%
_{i}\sum_{j=1}^{i+1}p_{j}a_{j}+P_{i}\sum_{j=i+1}^{n-1}p_{j+1}a_{j+1} \\
& =\left[ \left( P_{n}-P_{i}\right) P_{i+1}-P_{i}\left( P_{n}-P_{i+1}\right) %
\right] a_{i+1}+P_{i}\sum_{j=i+1}^{n-1}p_{j+1}a_{j+1}-\bar{P}%
_{i}\sum_{j=1}^{i+1}p_{j}a_{j} \\
& =P_{n}p_{i+1}a_{i+1}+P_{i}\sum_{j=i+1}^{n-1}p_{j+1}a_{j+1}-\bar{P}%
_{i}\sum_{j=1}^{i+1}p_{j}a_{j} \\
& =\left( P_{i}+\bar{P}_{i}\right)
p_{i+1}a_{i+1}+P_{i}\sum_{j=i+1}^{n-1}p_{j+1}a_{j+1}-\bar{P}%
_{i}\sum_{j=1}^{i+1}p_{j}a_{j} \\
& =P_{i}\sum_{j=i+1}^{n-1}p_{j}a_{j}-\bar{P}_{i}\sum_{j=1}^{i}p_{j}a_{j} \\
& =P_{i}\bar{A}_{i}\left( \mathbf{p}\right) -\bar{P}_{i}A_{i}\left( \mathbf{p%
}\right) \\
& =P_{i}A_{n}\left( \mathbf{p}\right) -P_{n}A_{i}\left( \mathbf{p}\right) ,
\end{align*}%
and the identity is proved.
\end{proof}

We are able now to state and prove the second identity for the \v{C}eby\v{s}%
ev functional.

\begin{theorem}
\label{t2.3}With the assumptions of Theorem \ref{t2.1}, we have the identity%
\begin{equation}
T_{n}\left( \mathbf{p};\mathbf{a},\mathbf{b}\right)
=\sum_{i=1}^{n-1}\sum_{j=1}^{n-1}P_{\min \left\{ i,j\right\} }\bar{P}_{\max
\left\{ i,j\right\} }\cdot \left\langle \Delta a_{j},\Delta
b_{i}\right\rangle .  \label{2.7}
\end{equation}
\end{theorem}

\begin{proof}
Follows by Theorem \ref{t2.1} and Lemma \ref{l2.2} and we omit the details.
\end{proof}

\section{New Inequalities}

The following result holds.

\begin{theorem}
\label{t3.1}Let $\left( H;\left\langle \cdot ,\cdot \right\rangle \right) $
be an inner product space over the real or complex number field $\mathbb{K}$%
; $\mathbf{p}=\left( p_{1},\dots ,p_{n}\right) \in \mathbb{R}^{n}$ and $%
\mathbf{a}=\left( a_{1},\dots ,a_{n}\right) ,\mathbf{b}=\left( b_{1},\dots
,b_{n}\right) \in H^{n}.$ Then we have the inequalities%
\begin{equation}
\left\vert T_{n}\left( \mathbf{p};\mathbf{a},\mathbf{b}\right) \right\vert
\leq \left\{ 
\begin{array}{l}
\max\limits_{1\leq i\leq n-1}\left\Vert P_{i}A_{n}\left( \mathbf{p}\right)
-P_{n}A_{i}\left( \mathbf{p}\right) \right\Vert
\sum\limits_{j=1}^{n-1}\left\Vert \Delta b_{j}\right\Vert ; \\ 
\\ 
\left( \sum\limits_{i=1}^{n-1}\left\Vert P_{i}A_{n}\left( \mathbf{p}\right)
-P_{n}A_{i}\left( \mathbf{p}\right) \right\Vert ^{q}\right) ^{\frac{1}{q}%
}\left( \sum\limits_{j=1}^{n-1}\left\Vert \Delta b_{j}\right\Vert
^{p}\right) ^{\frac{1}{p}} \\ 
\hfill \text{for \ }p>1,\ \frac{1}{p}+\frac{1}{q}=1; \\ 
\\ 
\sum\limits_{i=1}^{n-1}\left\Vert P_{i}A_{n}\left( \mathbf{p}\right)
-P_{n}A_{i}\left( \mathbf{p}\right) \right\Vert \cdot \max\limits_{1\leq
j\leq n-1}\left\Vert \Delta b_{j}\right\Vert .%
\end{array}%
\right.  \label{3.1}
\end{equation}%
All the inequalities in (\ref{3.1}) are sharp in the sense that the
constants $1$ cannot be replaced by smaller constants.
\end{theorem}

\begin{proof}
Using the first identity in (\ref{2.1}) and Schwarz's inequality in $H,$
i.e., $\left\vert \left\langle u,v\right\rangle \right\vert \leq \left\Vert
u\right\Vert \left\Vert v\right\Vert ,$ $u,v\in H,$ we have successively:%
\begin{align*}
\left\vert T_{n}\left( \mathbf{p};\mathbf{a},\mathbf{b}\right) \right\vert &
\leq \sum\limits_{i=1}^{n-1}\left\vert \left\langle P_{i}A_{n}\left( \mathbf{%
p}\right) -P_{n}A_{i}\left( \mathbf{p}\right) ,\Delta b_{i}\right\rangle
\right\vert \\
& \leq \sum\limits_{i=1}^{n-1}\left\Vert P_{i}A_{n}\left( \mathbf{p}\right)
-P_{n}A_{i}\left( \mathbf{p}\right) \right\Vert \left\Vert \Delta
b_{i}\right\Vert .
\end{align*}%
Using H\"{o}lder's inequality, we deduce the desired result (\ref{3.1}).

Let us prove, for instance, that the constant 1 in the second inequality is
best possible.

Assume, for $c>0,$ we have that%
\begin{equation}
\left\vert T_{n}\left( \mathbf{p};\mathbf{a},\mathbf{b}\right) \right\vert
\leq c\left( \sum\limits_{i=1}^{n-1}\left\Vert P_{i}A_{n}\left( \mathbf{p}%
\right) -P_{n}A_{i}\left( \mathbf{p}\right) \right\Vert ^{q}\right) ^{\frac{1%
}{q}}\left( \sum\limits_{j=1}^{n-1}\left\Vert \Delta b_{j}\right\Vert
^{p}\right) ^{\frac{1}{p}}  \label{3.2}
\end{equation}%
for $p>1,$ $\frac{1}{p}+\frac{1}{q}=1,$ $n\geq 2.$

If we choose $n=2,$ then we get%
\begin{equation*}
T_{2}\left( \mathbf{p};\mathbf{a},\mathbf{b}\right) =p_{1}p_{2}\left\langle
a_{2}-a_{1},b_{2}-b_{1}\right\rangle .
\end{equation*}%
Also, for $n=2,$%
\begin{equation*}
\left( \sum\limits_{i=1}^{n-1}\left\Vert P_{i}A_{n}\left( \mathbf{p}\right)
-P_{n}A_{i}\left( \mathbf{p}\right) \right\Vert ^{q}\right) ^{\frac{1}{q}%
}=\left\vert p_{1}p_{2}\right\vert \left\Vert a_{2}-a_{1}\right\Vert
\end{equation*}%
and 
\begin{equation*}
\left( \sum\limits_{j=1}^{n-1}\left\Vert \Delta b_{j}\right\Vert ^{p}\right)
^{\frac{1}{p}}=\left\Vert b_{2}-b_{1}\right\Vert ,
\end{equation*}%
and then, from (\ref{3.2}), for $n=2,$ we deduce%
\begin{equation}
\left\vert p_{1}p_{2}\right\vert \left\vert \left\langle
a_{2}-a_{1},b_{2}-b_{1}\right\rangle \right\vert \leq c\left\vert
p_{1}p_{2}\right\vert \left\Vert a_{2}-a_{1}\right\Vert \left\Vert
b_{2}-b_{1}\right\Vert .  \label{3.3}
\end{equation}%
If in (\ref{3.3}) we choose $a_{2}=b_{2},$ $a_{2}=b_{1}$ and $b_{2}\neq
b_{1},$ $p_{1},p_{2}\neq 0,$ we deduce $c\geq 1,$ proving that 1 is the best
possible constant in that inequality.
\end{proof}

The following corollary for the uniform distribution of the probability $%
\mathbf{p}$ holds.

\begin{corollary}
\label{c3.2}With the assumptions of Theorem \ref{t3.1} for $\mathbf{a}$ and $%
\mathbf{b},$ we have the inequalities%
\begin{equation}
0\leq \left\vert T_{n}\left( \mathbf{a},\mathbf{b}\right) \right\vert \leq 
\frac{1}{n^{2}}\times \left\{ 
\begin{array}{l}
\max\limits_{1\leq i\leq n-1}\left\Vert
i\sum\limits_{k=1}^{n}a_{k}-n\sum\limits_{k=1}^{i}a_{k}\right\Vert
\sum\limits_{j=1}^{n-1}\left\Vert \Delta b_{j}\right\Vert ; \\ 
\\ 
\left( \sum\limits_{i=1}^{n-1}\left\Vert
i\sum\limits_{k=1}^{n}a_{k}-n\sum\limits_{k=1}^{i}a_{k}\right\Vert
^{q}\right) ^{\frac{1}{q}}\left( \sum\limits_{j=1}^{n-1}\left\Vert \Delta
b_{j}\right\Vert ^{p}\right) ^{\frac{1}{p}} \\ 
\hfill \text{for \ }p>1,\ \frac{1}{p}+\frac{1}{q}=1; \\ 
\\ 
\sum\limits_{i=1}^{n-1}\left\Vert
i\sum\limits_{k=1}^{n}a_{k}-n\sum\limits_{k=1}^{i}a_{k}\right\Vert \cdot
\max\limits_{1\leq j\leq n-1}\left\Vert \Delta b_{j}\right\Vert .%
\end{array}%
\right.  \label{3.4}
\end{equation}
\end{corollary}

The following result may be stated as well.

\begin{theorem}
\label{t3.3}With the assumptions of Theorem \ref{t3.1} and if $P_{i}\neq 0$ $%
\left( i=1,\dots ,n\right) ,$ then we have the inequalities%
\begin{multline}
\left\vert T_{n}\left( \mathbf{p};\mathbf{a},\mathbf{b}\right) \right\vert 
\label{3.5} \\
\leq \left\vert P_{n}\right\vert \times \left\{ 
\begin{array}{l}
\max\limits_{1\leq i\leq n-1}\left\Vert \dfrac{1}{P_{n}}A_{n}\left( \mathbf{p%
}\right) -\dfrac{1}{P_{i}}A_{i}\left( \mathbf{p}\right) \right\Vert
\sum\limits_{i=1}^{n-1}\left\vert P_{i}\right\vert \left\Vert \Delta
b_{i}\right\Vert ; \\ 
\\ 
\left( \sum\limits_{i=1}^{n-1}\left\vert P_{i}\right\vert \left\Vert \dfrac{1%
}{P_{n}}A_{n}\left( \mathbf{p}\right) -\dfrac{1}{P_{i}}A_{i}\left( \mathbf{p}%
\right) \right\Vert ^{q}\right) ^{\frac{1}{q}}\left(
\sum\limits_{i=1}^{n-1}\left\vert P_{i}\right\vert \left\Vert \Delta
b_{i}\right\Vert ^{p}\right) ^{\frac{1}{p}} \\ 
\hfill \text{for \ }p>1,\ \frac{1}{p}+\frac{1}{q}=1; \\ 
\\ 
\sum\limits_{i=1}^{n-1}\left\vert P_{i}\right\vert \left\Vert \dfrac{1}{P_{n}%
}A_{n}\left( \mathbf{p}\right) -\dfrac{1}{P_{i}}A_{i}\left( \mathbf{p}%
\right) \right\Vert \cdot \max\limits_{1\leq i\leq n-1}\left\Vert \Delta
b_{i}\right\Vert .%
\end{array}%
\right. 
\end{multline}%
All the inequalities in (\ref{3.5}) are sharp in the sense that the constant
1 cannot be replaced by a smaller constant.
\end{theorem}

\begin{proof}
Using the second equality in (\ref{2.1}) and Schwarz's inequality, we have%
\begin{align*}
\left\vert T_{n}\left( \mathbf{p};\mathbf{a},\mathbf{b}\right) \right\vert &
\leq \left\vert P_{n}\right\vert \sum\limits_{i=1}^{n-1}\left\vert
P_{i}\right\vert \left\vert \left\langle \dfrac{1}{P_{n}}A_{n}\left( \mathbf{%
p}\right) -\dfrac{1}{P_{i}}A_{i}\left( \mathbf{p}\right) ,\Delta
b_{i}\right\rangle \right\vert \\
& \leq \left\vert P_{n}\right\vert \sum\limits_{i=1}^{n-1}\left\vert
P_{i}\right\vert \left\Vert \dfrac{1}{P_{n}}A_{n}\left( \mathbf{p}\right) -%
\dfrac{1}{P_{i}}A_{i}\left( \mathbf{p}\right) \right\Vert \left\Vert \Delta
b_{i}\right\Vert .
\end{align*}%
Using H\"{o}lder's weighted inequality, we deduce (\ref{3.5}).

The sharpness of the constant may be proven in a similar manner to the one
in Theorem \ref{t3.1}. We omit the details.
\end{proof}

The following corollary containing the unweighted inequalities holds.

\begin{corollary}
\label{c3.4}With the above assumptions for $\mathbf{a}$ and $\mathbf{b},$
one has%
\begin{equation}
\left\vert T_{n}\left( \mathbf{a},\mathbf{b}\right) \right\vert \leq \frac{1%
}{n}\times \left\{ 
\begin{array}{l}
\max\limits_{1\leq i\leq n-1}\left\Vert \dfrac{1}{n}\sum%
\limits_{k=1}^{n}a_{k}-\dfrac{1}{i}\sum\limits_{k=1}^{i}a_{k}\right\Vert
\sum\limits_{i=1}^{n-1}i\left\Vert \Delta b_{i}\right\Vert ; \\ 
\\ 
\left( \sum\limits_{i=1}^{n-1}i\left\Vert \dfrac{1}{n}\sum%
\limits_{k=1}^{n}a_{k}-\dfrac{1}{i}\sum\limits_{k=1}^{i}a_{k}\right\Vert
^{q}\right) ^{\frac{1}{q}}\left( \sum\limits_{i=1}^{n-1}i\left\Vert \Delta
b_{i}\right\Vert ^{p}\right) ^{\frac{1}{p}} \\ 
\hfill \text{for \ }p>1,\ \frac{1}{p}+\frac{1}{q}=1; \\ 
\\ 
\sum\limits_{i=1}^{n-1}i\left\Vert \dfrac{1}{n}\sum\limits_{k=1}^{n}a_{k}-%
\dfrac{1}{i}\sum\limits_{k=1}^{i}a_{k}\right\Vert \cdot \max\limits_{1\leq
i\leq n-1}\left\Vert \Delta b_{i}\right\Vert .%
\end{array}%
\right.  \label{3.6}
\end{equation}%
The inequalities (\ref{3.6}) are sharp in the sense mentioned above.
\end{corollary}

Another type of inequality may be stated if ones used the third identity in (%
\ref{2.1}) and H\"{o}lder's weighted inequality with the weights: $%
\left\vert P_{i}\right\vert \left\vert \bar{P}_{i}\right\vert ,$ $i\in
\left\{ 1,\dots ,n-1\right\} .$

\begin{theorem}
\label{t3.5}With the assumptions in Theorem \ref{t3.1} and if $P_{i},$ $\bar{%
P}_{i}\neq 0,$ $i\in \left\{ 1,\dots ,n-1\right\} ,$ then we have the
inequalities%
\begin{multline}
\left\vert T_{n}\left( \mathbf{p};\mathbf{a},\mathbf{b}\right) \right\vert
\label{3.7} \\
\leq \left\vert P_{n}\right\vert \times \left\{ 
\begin{array}{l}
\max\limits_{1\leq i\leq n-1}\left\Vert \dfrac{1}{\bar{P}_{i}}\bar{A}%
_{i}\left( \mathbf{p}\right) -\dfrac{1}{P_{i}}A_{i}\left( \mathbf{p}\right)
\right\Vert \sum\limits_{i=1}^{n-1}\left\vert P_{i}\right\vert \left\vert 
\bar{P}_{i}\right\vert \left\Vert \Delta b_{i}\right\Vert ; \\ 
\\ 
\left( \sum\limits_{i=1}^{n-1}\left\vert P_{i}\right\vert \left\vert \bar{P}%
_{i}\right\vert \left\Vert \dfrac{1}{\bar{P}_{i}}\bar{A}_{i}\left( \mathbf{p}%
\right) -\dfrac{1}{P_{i}}A_{i}\left( \mathbf{p}\right) \right\Vert
^{q}\right) ^{\frac{1}{q}}\left( \sum\limits_{i=1}^{n-1}\left\vert
P_{i}\right\vert \left\vert \bar{P}_{i}\right\vert \left\Vert \Delta
b_{i}\right\Vert ^{p}\right) ^{\frac{1}{p}} \\ 
\hfill \text{for \ }p>1,\ \frac{1}{p}+\frac{1}{q}=1; \\ 
\\ 
\sum\limits_{i=1}^{n-1}\left\vert P_{i}\right\vert \left\vert \bar{P}%
_{i}\right\vert \left\Vert \dfrac{1}{\bar{P}_{i}}\bar{A}_{i}\left( \mathbf{p}%
\right) -\dfrac{1}{P_{i}}A_{i}\left( \mathbf{p}\right) \right\Vert \cdot
\max\limits_{1\leq i\leq n-1}\left\Vert \Delta b_{i}\right\Vert .%
\end{array}%
\right.
\end{multline}%
In particular, if $p_{i}=\frac{1}{n},$ $i\in \left\{ 1,\dots ,n\right\} ,$
then we have%
\begin{multline}
\left\vert T_{n}\left( \mathbf{a},\mathbf{b}\right) \right\vert  \label{3.8}
\\
\leq \frac{1}{n^{2}}\times \left\{ 
\begin{array}{l}
\max\limits_{1\leq i\leq n-1}\left\Vert \dfrac{1}{n-i}\sum%
\limits_{k=i+1}^{n}a_{k}-\dfrac{1}{i}\sum\limits_{k=1}^{i}a_{k}\right\Vert
\sum\limits_{i=1}^{n-1}i\left( n-i\right) \left\Vert \Delta b_{i}\right\Vert
; \\ 
\\ 
\left( \sum\limits_{i=1}^{n-1}i\left( n-i\right) \left\Vert \dfrac{1}{n-i}%
\sum\limits_{k=i+1}^{n}a_{k}-\dfrac{1}{i}\sum\limits_{k=1}^{i}a_{k}\right%
\Vert ^{q}\right) ^{\frac{1}{q}}\left( \sum\limits_{i=1}^{n-1}i\left(
n-i\right) \left\Vert \Delta b_{i}\right\Vert ^{p}\right) ^{\frac{1}{p}} \\ 
\hfill \text{for \ }p>1,\ \frac{1}{p}+\frac{1}{q}=1; \\ 
\\ 
\sum\limits_{i=1}^{n-1}i\left( n-i\right) \left\Vert \dfrac{1}{n-i}%
\sum\limits_{k=i+1}^{n}a_{k}-\dfrac{1}{i}\sum\limits_{k=1}^{i}a_{k}\right%
\Vert \cdot \max\limits_{1\leq i\leq n-1}\left\Vert \Delta b_{i}\right\Vert .%
\end{array}%
\right.
\end{multline}%
The inequalities in (\ref{3.7}) and (\ref{3.8}) are sharp in the above
mentioned sense.
\end{theorem}

A different approach may be considered if one uses the representation in
terms of double sums for the \v{C}eby\v{s}ev functional provided by Theorem %
\ref{t2.3}.

The following result holds.

\begin{theorem}
\label{t3.6}With the above assumptions of Theorem \ref{t3.1}, we have the
inequalities%
\begin{multline}
\left\vert T_{n}\left( \mathbf{p};\mathbf{a},\mathbf{b}\right) \right\vert
\label{3.9} \\
\leq \left\vert P_{n}\right\vert \times \left\{ 
\begin{array}{l}
\max\limits_{1\leq i,j\leq n-1}\left\{ \left\vert P_{\min \left\{
i,j\right\} }\right\vert ,\left\vert \bar{P}_{\max \left\{ i,j\right\}
}\right\vert \right\} \sum\limits_{i=1}^{n-1}\left\Vert \Delta
a_{i}\right\Vert \sum\limits_{i=1}^{n-1}\left\Vert \Delta b_{i}\right\Vert ;
\\ 
\\ 
\left( \sum\limits_{i=1}^{n-1}\sum\limits_{j=1}^{n-1}\left\vert P_{\min
\left\{ i,j\right\} }\right\vert ^{q}\left\vert \bar{P}_{\max \left\{
i,j\right\} }\right\vert ^{q}\right) ^{\frac{1}{q}}\left(
\sum\limits_{i=1}^{n-1}\left\Vert \Delta a_{i}\right\Vert ^{p}\right) ^{%
\frac{1}{p}}\left( \sum\limits_{i=1}^{n-1}\left\Vert \Delta b_{i}\right\Vert
^{p}\right) ^{\frac{1}{p}} \\ 
\hfill \text{for \ }p>1,\ \frac{1}{p}+\frac{1}{q}=1; \\ 
\\ 
\sum\limits_{i=1}^{n-1}\sum\limits_{j=1}^{n-1}\left\vert P_{\min \left\{
i,j\right\} }\right\vert \left\vert \bar{P}_{\max \left\{ i,j\right\}
}\right\vert \max\limits_{1\leq i\leq n-1}\left\Vert \Delta a_{i}\right\Vert
\max\limits_{1\leq i\leq n-1}\left\Vert \Delta b_{i}\right\Vert .%
\end{array}%
\right.
\end{multline}%
The inequalities are sharp in the sense mentioned above.
\end{theorem}

The proof follows by the identity (\ref{2.7}) on using H\"{o}lder's
inequality for double sums and we omit the details.

Now, define%
\begin{equation*}
k_{\infty }:=\max_{1\leq i,j\leq n-1}\left\{ \frac{\min \left\{ i,j\right\} 
}{n}\left( 1-\frac{\max \left\{ i,j\right\} }{n}\right) \right\} ,\ \ n\geq
2.
\end{equation*}%
Using the elementary inequality%
\begin{equation*}
ab\leq \frac{1}{4}\left( a+b\right) ^{2},\ \ \ \ a,b\in \mathbb{R};
\end{equation*}%
we deduce%
\begin{equation*}
\min \left\{ i,j\right\} \left( n-\max \left\{ i,j\right\} \right) \leq 
\frac{1}{4}\left( n-\left\vert i-j\right\vert \right) ^{2}
\end{equation*}%
for $1\leq i,j\leq n-1.$ Consequently, we have%
\begin{equation*}
k_{\infty }\leq \frac{1}{4n^{2}}\max_{1\leq i,j\leq n-1}\left\{ \left(
n-\left\vert i-j\right\vert \right) ^{2}\right\} =\frac{1}{4}.
\end{equation*}

We may now state the following corollary of Theorem \ref{t3.6}.

\begin{corollary}
\label{c3.7}With the assumptions of Theorem \ref{t3.1} for $\mathbf{a}$ and $%
\mathbf{b},$ we have the inequality%
\begin{align}
\left\vert T_{n}\left( \mathbf{a},\mathbf{b}\right) \right\vert & \leq
k_{\infty }\sum\limits_{i=1}^{n-1}\left\Vert \Delta a_{i}\right\Vert
\sum\limits_{i=1}^{n-1}\left\Vert \Delta b_{i}\right\Vert  \label{3.10} \\
& \leq \frac{1}{4}\sum\limits_{i=1}^{n-1}\left\Vert \Delta a_{i}\right\Vert
\sum\limits_{i=1}^{n-1}\left\Vert \Delta b_{i}\right\Vert .  \notag
\end{align}%
The constant $\frac{1}{4}$ cannot be replaced in general by a smaller
constant.
\end{corollary}

\begin{remark}
\label{r3.8}The inequality (\ref{3.10}) is better than the third inequality
in (\ref{1.10}).
\end{remark}

Consider now, for $q>1,$ the number%
\begin{equation*}
k_{q}:=\frac{1}{n^{2}}\left( \sum\limits_{i=1}^{n-1}\sum\limits_{j=1}^{n-1}%
\left[ \min \left\{ i,j\right\} \left( n-\max \left\{ i,j\right\} \right) %
\right] ^{q}\right) ^{\frac{1}{q}}.
\end{equation*}%
We observe, by symmetry of the terms under the summation symbol, we have that%
\begin{equation*}
k_{q}=\frac{1}{n^{2}}\left( 2\sum_{1\leq i<j\leq n-1}i^{q}\left( n-j\right)
^{q}+\sum_{i=1}^{n-1}i^{q}\left( n-i\right) ^{q}\right) ^{\frac{1}{q}},
\end{equation*}%
that may be computed exactly if $q=2$ or another natural number.

Since, as above,%
\begin{equation*}
\left[ \min \left\{ i,j\right\} \left( n-\max \left\{ i,j\right\} \right) %
\right] ^{q}\leq \frac{1}{4^{q}}\left( n-\left\vert i-j\right\vert \right)
^{2q},
\end{equation*}%
we deduce%
\begin{align*}
k_{q}& \leq \frac{1}{4n^{2}}\left(
\sum\limits_{i=1}^{n-1}\sum\limits_{j=1}^{n-1}\left( n-\left\vert
i-j\right\vert \right) ^{2q}\right) ^{\frac{1}{q}} \\
& \leq \frac{1}{4n^{2}}\left[ \left( n-1\right) ^{2}n^{2q}\right] ^{\frac{1}{%
q}} \\
& =\frac{1}{4}\left( n-1\right) ^{\frac{2}{q}}.
\end{align*}%
Consequently, we may state the following corollary as well.

\begin{corollary}
\label{c3.9}With the assumptions of Theorem \ref{t3.1} for $\mathbf{a}$ and $%
\mathbf{b},$ we have the inequalities%
\begin{align}
\left\vert T_{n}\left( \mathbf{a},\mathbf{b}\right) \right\vert & \leq
k_{q}\left( \sum\limits_{i=1}^{n-1}\left\Vert \Delta a_{i}\right\Vert
^{p}\right) ^{\frac{1}{p}}\left( \sum\limits_{i=1}^{n-1}\left\Vert \Delta
b_{i}\right\Vert ^{p}\right) ^{\frac{1}{p}}  \label{3.11} \\
& \leq \frac{1}{4}\left( n-1\right) ^{\frac{2}{q}}\left(
\sum\limits_{i=1}^{n-1}\left\Vert \Delta a_{i}\right\Vert ^{p}\right) ^{%
\frac{1}{p}}\left( \sum\limits_{i=1}^{n-1}\left\Vert \Delta b_{i}\right\Vert
^{p}\right) ^{\frac{1}{p}},  \notag
\end{align}%
provided $p>1$ , $\frac{1}{p}+\frac{1}{q}=1.$ The constant $\frac{1}{4}$
cannot be replaced in general by a smaller constant.
\end{corollary}

Finally, if we denote%
\begin{equation*}
k_{1}:=\frac{1}{n^{2}}\sum\limits_{i=1}^{n-1}\sum\limits_{j=1}^{n-1}\left[
\min \left\{ i,j\right\} \left( n-\max \left\{ i,j\right\} \right) \right] ,
\end{equation*}%
then we observe, for $\mathbf{u}=\left( \frac{1}{n},\dots ,\frac{1}{n}%
\right) ,$ $\mathbf{e}=\left( 1,2,\dots ,n\right) ,$ that%
\begin{equation*}
k_{1}=T_{n}\left( \mathbf{u};\mathbf{e},\mathbf{e}\right) =\frac{1}{n}%
\sum\limits_{i=1}^{n}i^{2}-\left( \frac{1}{n}\sum\limits_{i=1}^{n}i\right)
^{2}=\frac{1}{12}\left( n^{2}-1\right) ,
\end{equation*}%
and by Theorem \ref{t3.6}, we deduce the inequality%
\begin{equation*}
\left\vert T_{n}\left( \mathbf{a},\mathbf{b}\right) \right\vert \leq \frac{1%
}{12}\left( n^{2}-1\right) \max_{1\leq j\leq n-1}\left\Vert \Delta
a_{j}\right\Vert \max_{1\leq j\leq n-1}\left\Vert \Delta b_{j}\right\Vert .
\end{equation*}%
Note that, the above inequality has been discovered using a different method
in \cite{SSD3}. The constant $\frac{1}{12}$ is best possible.

\section{Reverses for Jensen's Inequality}

Let $\left( H;\left\langle \cdot ,\cdot \right\rangle \right) $ be a real
inner product space and $F:H\rightarrow \mathbb{R}$ a Fr\'{e}chet
differentiable convex function on $H.$ If $\triangledown F:H\rightarrow H$
denotes the gradient operator associated to $F,$ then we have the inequality%
\begin{equation*}
F\left( x\right) -F\left( y\right) \geq \left\langle \triangledown F\left(
y\right) ,x-y\right\rangle
\end{equation*}%
for each $x,y\in H.$

The following result has been obtained in \cite{SSD3}.

\begin{theorem}
\label{t4.1}Let $F:H\rightarrow \mathbb{R}$ be as above and $z_{i}\in H,$ $%
i\in \left\{ 1,\dots ,n\right\} .$ If $q_{i}\geq 0$ \ $\left( i\in \left\{
1,\dots ,n\right\} \right) $ with $\sum_{i=1}^{n}q_{i}=1,$ then we have the
following reverse of Jensen's inequality%
\begin{align}
0& \leq \sum_{i=1}^{n}q_{i}F\left( z_{i}\right) -F\left(
\sum_{i=1}^{n}q_{i}z_{i}\right)  \label{4.2} \\
& \leq \left\{ 
\begin{array}{l}
\left[ \sum\limits_{i=1}^{n}i^{2}q_{i}-\left(
\sum\limits_{i=1}^{n}iq_{i}\right) ^{2}\right] \max\limits_{k=1,\dots
,n-1}\left\Vert \Delta \left( \triangledown F\left( z_{i}\right) \right)
\right\Vert \max\limits_{k=1,\dots ,n-1}\left\Vert \Delta z_{i}\right\Vert ;
\\ 
\\ 
\left[ \sum\limits_{1\leq j<i\leq n}q_{i}q_{j}\left( i-j\right) \right]
\left( \sum\limits_{i=1}^{n-1}\left\Vert \Delta \left( \triangledown F\left(
z_{i}\right) \right) \right\Vert ^{p}\right) ^{\frac{1}{p}}\left(
\sum\limits_{i=1}^{n-1}\left\Vert \Delta z_{i}\right\Vert ^{q}\right) ^{%
\frac{1}{q}} \\ 
\hfill \text{if }p>1,\ \frac{1}{p}+\frac{1}{q}=1; \\ 
\\ 
\dfrac{1}{2}\left[ \sum\limits_{i=1}^{n}q_{i}\left( 1-q_{i}\right) \right]
\sum\limits_{i=1}^{n-1}\left\Vert \Delta \left( \triangledown F\left(
z_{i}\right) \right) \right\Vert \sum\limits_{i=1}^{n-1}\left\Vert \Delta
z_{i}\right\Vert .%
\end{array}%
\right.  \notag \\
&  \notag
\end{align}
\end{theorem}

The unweighted case may useful in application and is incorporated in the
following corollary.

\begin{corollary}
Let $F:H\rightarrow \mathbb{R}$ be as above and $z_{i}\in H,$ $i\in \left\{
1,\dots ,n\right\} .$ Then we have the inequalities%
\begin{align*}
0& \leq \frac{1}{n}\sum_{i=1}^{n}F\left( z_{i}\right) -F\left( \frac{1}{n}%
\sum_{i=1}^{n}z_{i}\right) \\
& \leq \left\{ 
\begin{array}{l}
\dfrac{n^{2}-1}{12}\max\limits_{k=1,\dots ,n-1}\left\Vert \Delta \left(
\triangledown F\left( z_{k}\right) \right) \right\Vert
\max\limits_{k=1,\dots ,n-1}\left\Vert \Delta z_{k}\right\Vert ; \\ 
\\ 
\dfrac{n^{2}-1}{6n}\left( \sum\limits_{k=1}^{n-1}\left\Vert \Delta \left(
\triangledown F\left( z_{k}\right) \right) \right\Vert ^{p}\right) ^{\frac{1%
}{p}}\left( \sum\limits_{k=1}^{n-1}\left\Vert \Delta z_{k}\right\Vert
^{q}\right) ^{\frac{1}{q}} \\ 
\hfill \text{if }p>1,\ \frac{1}{p}+\frac{1}{q}=1; \\ 
\\ 
\dfrac{n-1}{2n}\sum\limits_{k=1}^{n-1}\left\Vert \Delta \left( \triangledown
F\left( z_{k}\right) \right) \right\Vert \sum\limits_{k=1}^{n-1}\left\Vert
\Delta z_{k}\right\Vert .%
\end{array}%
\right. \\
&
\end{align*}
\end{corollary}

By making use, of Theorem \ref{t3.6}, we can state the following result as
well:

\begin{theorem}
\label{t4.2}Let $F:H\rightarrow \mathbb{R}$ be as above and $z_{i}\in H,$ $%
i\in \left\{ 1,\dots ,n\right\} .$ If $q_{i}\geq 0$ \ $\left( i\in \left\{
1,\dots ,n\right\} \right) $ with $\sum_{i=1}^{n}q_{i}=1,$ then we have the
following reverse of Jensen's inequality%
\begin{align}
0& \leq \sum_{i=1}^{n}q_{i}F\left( z_{i}\right) -F\left(
\sum_{i=1}^{n}q_{i}z_{i}\right)  \notag \\
&  \label{4.3} \\
& \leq \left\{ 
\begin{array}{l}
\max\limits_{1\leq i,j\leq n-1}\left\{ Q_{\min \left\{ i,j\right\} },%
\overline{Q}_{\max \left\{ i,j\right\} }\right\}
\sum\limits_{i=1}^{n-1}\left\Vert \Delta \left( \triangledown F\left(
z_{i}\right) \right) \right\Vert \sum\limits_{i=1}^{n-1}\left\Vert \Delta
z_{i}\right\Vert ; \\ 
\\ 
\left( \sum\limits_{i=1}^{n-1}\sum\limits_{j=1}^{n-1}Q_{\min \left\{
i,j\right\} }^{q}\overline{Q}_{\max \left\{ i,j\right\} }^{q}\right) ^{\frac{%
1}{q}}\left( \sum\limits_{i=1}^{n-1}\left\Vert \Delta \left( \triangledown
F\left( z_{i}\right) \right) \right\Vert ^{p}\right) ^{\frac{1}{p}}\left(
\sum\limits_{i=1}^{n-1}\left\Vert \Delta z_{i}\right\Vert ^{p}\right) ^{%
\frac{1}{p}} \\ 
\hfill \text{for \ }p>1,\ \frac{1}{p}+\frac{1}{q}=1; \\ 
\\ 
\sum\limits_{i=1}^{n-1}\sum\limits_{j=1}^{n-1}Q_{\min \left\{ i,j\right\} }%
\overline{Q}_{\max \left\{ i,j\right\} }\max\limits_{1\leq i\leq
n-1}\left\Vert \Delta \left( \triangledown F\left( z_{i}\right) \right)
\right\Vert \max\limits_{1\leq i\leq n-1}\left\Vert \Delta z_{i}\right\Vert .%
\end{array}%
\right.  \notag
\end{align}
\end{theorem}

\begin{proof}
We know, see for example \cite[Eq. (4.4)]{SSD1}, that the following reverse
of Jensen's inequality for Fr\'{e}chet differentiable convex functions%
\begin{align}
0& \leq \sum_{i=1}^{n}q_{i}F\left( z_{i}\right) -F\left(
\sum_{i=1}^{n}q_{i}z_{i}\right)  \label{4.4} \\
& \leq \sum_{i=1}^{n}q_{i}\left\langle \triangledown F\left( z_{i}\right)
,z_{i}\right\rangle -\left\langle \sum_{i=1}^{n}q_{i}\triangledown F\left(
z_{i}\right) ,\sum_{i=1}^{n}q_{i}z_{i}\right\rangle  \notag
\end{align}%
holds.

Now, if we apply Theorem \ref{t3.6} for the choices $a_{i}=\triangledown
F\left( z_{i}\right) ,b_{i}=z_{i}$ and $p_{i}=q_{i}\left( i=1,...,n\right) ,$
then we may state%
\begin{align}
& \sum_{i=1}^{n}q_{i}\left\langle \triangledown F\left( z_{i}\right)
,z_{i}\right\rangle -\left\langle \sum_{i=1}^{n}q_{i}\triangledown F\left(
z_{i}\right) ,\sum_{i=1}^{n}q_{i}z_{i}\right\rangle  \notag \\
&  \label{4.5} \\
& \leq \left\{ 
\begin{array}{l}
\max\limits_{1\leq i,j\leq n-1}\left\{ Q_{\min \left\{ i,j\right\} },%
\overline{Q}_{\max \left\{ i,j\right\} }\right\}
\sum\limits_{i=1}^{n-1}\left\Vert \Delta \left( \triangledown F\left(
z_{i}\right) \right) \right\Vert \sum\limits_{i=1}^{n-1}\left\Vert \Delta
z_{i}\right\Vert ; \\ 
\\ 
\left( \sum\limits_{i=1}^{n-1}\sum\limits_{j=1}^{n-1}Q_{\min \left\{
i,j\right\} }^{q}\overline{Q}_{\max \left\{ i,j\right\} }^{q}\right) ^{\frac{%
1}{q}}\left( \sum\limits_{i=1}^{n-1}\left\Vert \Delta \left( \triangledown
F\left( z_{i}\right) \right) \right\Vert ^{p}\right) ^{\frac{1}{p}}\left(
\sum\limits_{i=1}^{n-1}\left\Vert \Delta z_{i}\right\Vert ^{p}\right) ^{%
\frac{1}{p}} \\ 
\text{for \ }p>1,\ \frac{1}{p}+\frac{1}{q}=1; \\ 
\\ 
\sum\limits_{i=1}^{n-1}\sum\limits_{j=1}^{n-1}Q_{\min \left\{ i,j\right\} }%
\overline{Q}_{\max \left\{ i,j\right\} }\max\limits_{1\leq i\leq
n-1}\left\Vert \Delta \left( \triangledown F\left( z_{i}\right) \right)
\right\Vert \max\limits_{1\leq i\leq n-1}\left\Vert \Delta z_{i}\right\Vert .%
\end{array}%
\right.  \notag
\end{align}%
Finally, on making use of the inequalities (\ref{4.4}) and (\ref{4.5}), we
deduce the desired result (\ref{4.3}).
\end{proof}

The unweighted case may be useful in application and is incorporated in the
following corollary.

\begin{corollary}
Let $F:H\rightarrow \mathbb{R}$ be as above and $z_{i}\in H,$ $i\in \left\{
1,\dots ,n\right\} .$ Then we have the inequalities%
\begin{align*}
0& \leq \frac{1}{n}\sum_{i=1}^{n}F\left( z_{i}\right) -F\left( \frac{1}{n}%
\sum_{i=1}^{n}z_{i}\right) \\
& \leq \left\{ 
\begin{array}{l}
\dfrac{n^{2}-1}{12}\max\limits_{k=1,\dots ,n-1}\left\Vert \Delta \left(
\triangledown F\left( z_{k}\right) \right) \right\Vert
\max\limits_{k=1,\dots ,n-1}\left\Vert \Delta z_{k}\right\Vert ; \\ 
\\ 
\frac{1}{4}\left( n-1\right) ^{\frac{2}{q}}\left(
\sum\limits_{k=1}^{n-1}\left\Vert \Delta \left( \triangledown F\left(
z_{k}\right) \right) \right\Vert ^{p}\right) ^{\frac{1}{p}}\left(
\sum\limits_{k=1}^{n-1}\left\Vert \Delta z_{k}\right\Vert ^{p}\right) ^{%
\frac{1}{p}} \\ 
\hfill \text{if }p>1,\ \frac{1}{p}+\frac{1}{q}=1; \\ 
\\ 
\frac{1}{4}\sum\limits_{k=1}^{n-1}\left\Vert \Delta \left( \triangledown
F\left( z_{k}\right) \right) \right\Vert \sum\limits_{k=1}^{n-1}\left\Vert
\Delta z_{k}\right\Vert .%
\end{array}%
\right.
\end{align*}
\end{corollary}

\begin{remark}
If one applies the other Gr\"{u}ss' type inequalities obtained in the
previous section, then one can obtain other reverses for Jensen's discrete
inequality for convex functions defined on inner product spaces. We do not
present them here.
\end{remark}

\end{document}